\documentclass[11pt]{article}
\usepackage[T1]{fontenc}
\usepackage{latexsym,amssymb,amsmath,amsfonts,amsthm}
\usepackage{graphicx}
\usepackage{subfigure}
\usepackage{placeins}
\usepackage{color}
\usepackage{indentfirst}
\numberwithin{equation}{section}
\newcommand{\R}{{\mathbb R}}

\newcommand{\N}{{\mathbb N}}
\newcommand{\C}{{\mathbb C}}
\newcommand{\s}{{\mathbb S}}

\newcommand{\be}{\begin{eqnarray}}
\newcommand{\ben}{\begin{eqnarray*}}
\newcommand{\en}{\end{eqnarray}}
\newcommand{\enn}{\end{eqnarray*}}

\newtheorem{theorem}{Theorem}[section]
\newtheorem{lemma}[theorem]{Lemma}
\newtheorem{corollary}[theorem]{Corollary}
\newtheorem{definition}[theorem]{Definition}

\definecolor{rot}{rgb}{1.000,0.000,0.000}

\begin{document}
\renewcommand{\theequation}{\arabic{section}.\arabic{equation}}
\begin{titlepage}
\title{\bf Weakly singular corners always scatter}


\author{
Long Li\thanks{LMAM, School of Mathematical Sciences, Peking University, Beijing 100871,China, and
Beijing Advanced Innovation Center for Imaging Technology, Capital Normal University, Beijing 100048, China ({\tt 1601110046@pku.edu.cn}).
}
\and Guanghui Hu\thanks{Corresponding author: Beijing Computational Sciences Research Center, Beijing 100193, China ({\tt
hu@csrc.ac.cn}).}
\and Jiansheng Yang \thanks{LMAM, School of Mathematical Sciences, Peking University, Beijing 100871,China, and
Beijing Advanced Innovation Center for Imaging Technology, Capital Normal University, Beijing 100048, China ({\tt jsyang@pku.edu.cn}).
}
 }
\date{}
\end{titlepage}
\maketitle
\vspace{.2in}

\begin{abstract} Assume that a bounded
scatterer is embedded into an infinite homogeneous isotropic background medium in two dimensions. The refractive index function is supposed to be piecewise constant. If the scattering interface contains a weakly or strongly singular point, we prove that the scattered field cannot vanish identically. This particularly leads to the absence of non-scattering energies for piecewise analytic interfaces with a weakly singular point. Local uniqueness is obtained for shape identification problems in inverse medium scattering with a single far-field pattern.

\vspace{.2in} {\bf Keywords}: Uniqueness, inverse medium scattering, non-scattering energy; weakly singular corners.
\end{abstract}

\section{Introduction}
Assume a time-harmonic incoming wave $u^{in}$ is incident onto a bounded penetrable
scatterer $D\subset \R^2$ embedded in a homogeneous isotropic background medium.
 We assume that the boundary $\partial D$ is Lipschitz continuous and piecewise analytic, and that the complement $D^e:=\R^2\backslash\overline{D}$ of $D$ is connected.
 The wave propagation of the total field $u=u^{in}+u^{sc}$ is then modeled by the Helmholtz equation
\begin{equation}\label{eq:Helm}
\Delta u+k^2 q\, u=0\quad\mbox{in}\quad \R^2.
\end{equation}
In this paper the  refractive index (potential) function $q$ is supposed to be a piecewise constant function, given by
\ben
q(x)=\left\{\begin{array}{lll}
1,\quad\mbox{if}\quad x\in D^e,\\
q_0\neq 1,\quad\mbox{if}\quad x\in \overline{D}.
\end{array}\right.
\enn
Across the interface $\partial D$, we assume the continuity of the total field and its normal derivative, i.e.,
\be\label{TE}
u^+=u^-,\quad \partial_\nu u^+=\partial_\nu u^-\quad\mbox{on}\;\partial D.
\en
Here the superscripts $(\cdot)^\pm$ stand for the limits taken from outside and inside, respectively, and $\nu\in \s:=\{x\in\R^2: |x|=1\}$ is the unit normal on $\partial D$ pointing into $D^e$.
 At the infinity, the perturbed scattered field
 $u^{sc}$ is supposed to fulfill the Sommerfeld radiation condition
\begin{equation}\label{eq:radiation}
\lim_{|x|\rightarrow \infty} |x|\,\left\{ \frac{\partial u^{sc}}{\partial |x|}-ik u^{sc} \right\}=0.
\end{equation}
The unique solvability of the scattering problem \eqref{eq:Helm}, \eqref{eq:radiation} and (\ref{TE}) in $H^2_{loc}(\R^2)$ is well known (see e.g., \cite[Chapter 8]{CK}). In particular, the Sommerfeld radiation condition (\ref{eq:radiation}) leads to  the asymptotic expansion
\begin{equation}\label{eq:farfield}
u^{sc}(x)=\frac{e^{ik |x|}}{\sqrt{|x|}}\; u^\infty(\hat x)+\mathcal{O}\left(\frac{1}{|x|^{3/2}}\right),\quad |x|\rightarrow+\infty,
\end{equation}
 uniformly in all directions $\hat x:=x/|x|$, $x\in\mathbb{R}$. The function $u^\infty(\hat x)$ is an analytic function defined on $\s^{2}$ and is  referred to as the \emph{far-field pattern} or the \emph{scattering amplitude}.  The vector $\hat{x}\in\s$ is called the observation direction of the far field. The classical inverse scattering problem consists of the recovery of the boundary $\partial D$ from the far-field patterns corresponding to one or several incident plane waves. In this paper we are
 concerned with the following questions:
 \begin{description}
\item[(i)] Does the obstacle $D$ scatter any incident wave trivially (that is, $u^{sc}\equiv 0$) ?
\item[(ii)] Does the far-field pattern of a single incoming wave uniquely determine $\partial D$ ?
\end{description}
A negative answer to the first question means that acoustic cloaking cannot be achieved using isotropic materials, while a positive answer implies that $k^2$ is a non-scattering wavenumber (energy).
The study of non-scattering energies dates back to \cite{KS} in the case of a convex (planar) corner domain, where notion of \emph{scattering support} for an inhomogeneous medium was explored. In one of the authors' previous work \cite{EH2017}, it was shown that variable potential functions with the following corners on $\partial D$:
\begin{itemize}
\item curvilinear polygonal corners in $\R^2$;
\item curvilinear polyhedral corners in $\R^3$;
\item circular conic corners in $\R^3$;
\end{itemize}
scatters every incident wave non-trivially.  Earlier publications were devoted to the
absence of non-scattering energies under more restrictive assumptions on the smoothness of the potential or the angle of the corner. Here we mention the following works in the acoustic case:
\begin{itemize}
\item $C^\infty$-potentials with rectangular corners in $\R^n$ ($n\geq 2$) \cite{BLS};
\item H\"older continuous potentials with convex corners in $\R^2$, and with circular conic corners in $\R^3$ whose opening angle is outside of a countable subset of $(0,\pi)$ \cite{PSV};
\item analytical potentials with arbitrary polygonal corners or polyhedral wedge corners \cite{ElHu2015};
\item H\"older continuous potentials with rectangular corners in $\R^3$ \cite{HSV}.
\end{itemize}
 The argument of the pioneering work  \cite{BLS} was based on the use of complex geometric optics (CGO) solutions, which
 was later extended to \cite{PSV} and \cite{HSV} for treating less regular potentials and convex corners.
 The approach of \cite{ElHu2015} relies on the expansion of solutions to the Helmholtz equation with real-analytic potentials. For general potentials and corners, the absence of non-scattering energies can be verified via singularity analysis of the inhomogeneous Laplace equation in a cone \cite{EH2017}. We remark that the first question is closely related to the second one, that is, the approach for proving absence of non-scattering wavenumbers implies uniqueness to shape identification problems in inverse medium scattering.
 It was first proved in \cite{ElHu2015} that the shape of a convex penetrable obstacle of polygonal or polyhedral type with an unknown analytical potential can be uniquely determined by a single far-field pattern. The CGO-solution methods of \cite{PSV,BLS} also lead to uniqueness in shape identification but are confined so far to convex polygons in $\R^2$ and rectangular boxes in $\R^3$ with H\"older continuous potentials (see \cite{HSV}).
 In \cite[Corollay 2.1]{EH2017}, the uniqueness result of \cite{ElHu2015} was extended to
  more general potential functions using the data of a single far-field pattern.

\section{Main results}
The main purpose of this paper is to exclude (positive) real non-scattering energies when $\partial D$ contains a weakly singular corner, around which the boundary is allowed to be $C^1$-smooth but piecewise analytic.
The corners mentioned in the previous section are all strongly singular in the following sense.
\begin{definition}\label{def:SS}
A point $O\in \partial D\subset \R^2$ is called strongly singular if the boundary around $O$ can be locally parameterized by a continuous and piecewise analytic function whose derivative is discontinuous at $O$.
\end{definition}
Evidently,  every planar corner point of a polygon with flat slides is strongly singular, because the boundary can be locally parameterized by a piecewise linear function, whose
 first derivative is piecewise constant. A curvilinear corner of $D$ (see e.g., Definition 2.1 of \cite{EH2017} for a precise description) is also strongly singular by Definition \ref{def:SS}.  Below we state the definition of weakly singular corners to be explored within the scope of this paper.

\begin{definition}\label{def:WS}
A point $O=(0,0)\in \partial D$ is called weakly singular if the subboundary $B_\epsilon(O)\cap \partial D$ for some $\epsilon>0 $ can be parameterized by the  polynomial function $x_2=f(x_1)$, $x_1\in(-1,1)$, where
\be\label{f1}
f(x_1)=\left\{\begin{array}{lll}
c_1\, x_1^{\alpha_1}, &&1>x_1\geq 0,\\
c_2\, x_1^{\alpha_2},&& -1<x_1\leq 0,
\end{array}\right.
\en and the coefficients $c_j\in \R$ and $\alpha_j\in \N^+$ are assumed to fulfill the relations
\ben
(c_1, \alpha_1)\neq (c_2, \alpha_2),\;c_1^2+c_2^2\ne0,\;\;\alpha_j\geq 2.
\enn
The order of the singularity at $O$ is defined as
\ben
\beta:=\left\{\begin{array}{lll}
\min\{\alpha_1, \alpha_2\}&&\mbox{if}\quad c_1\neq 0, c_2\neq 0,\\
\alpha_1 &&\mbox{if}\quad c_2= 0,\\
\alpha_2 && \mbox{if}\quad c_1=0.
\end{array}\right.
\enn
\end{definition}
The boundary around  a weakly singular corner of order $\beta$ is  $C^{\beta-1}$-smooth but  piecewise  $C^{\beta}$-smooth, that is, the $\beta$-th derivative is discontinuous at O.
A singular point of order one  must be strongly singular in the sense of Definition \ref{def:SS}.
\begin{definition}\label{def:S}
A point $O\in \partial D$ is called singular if it is either strongly singular in the sense of Definition \ref{def:SS} or weakly singular in the sense of Definition \ref{def:WS}.
\end{definition}
The singular points defined by Definition \ref{def:S} form only a subset of non-analytic points of the boundary. In fact, the polynomial functions described in (\ref{f1}) can be regarded as the leading terms of the Taylor expansion of an analytic function at $x_1=0^\pm$.
If $q$ is a piecewise constant function in $\R^2$, we shall prove that
\begin{theorem}\label{TH}
The obstacle $D\subset \R^2$ scatters every incoming wave, if $\partial D$ contains a singular point.
\end{theorem}
 Note that when $\partial D$ possesses a strongly singular corner, Theorem \ref{TH} has been implicitly contained in \cite{EH2017}. The main contribution of this paper is to verify
Theorem \ref{TH} for weakly singular corners in $\R^2$. The above theorem implies that a Lipschitz domain with a singular point on the boundary scatters every incoming wave trivially in two dimensions. Theorem \ref{TH} follows straightforwardly from  Lemma \ref{lem:SS} for strongly singular corners and Lemma \ref{lem:WS}  for weakly singular corners.

Only local properties of the Helmholtz equation are involved in the proof of Theorem \ref{TH}.
Consequently, we get a local uniqueness result to the inverse scattering for shape identification:
 \begin{theorem}\label{TH2} Let $D_j$ ($j=1,2$) be two penetrable obstacles in $\R^2$ with the piecewise constant potential functions $q_j$, respectively.  If $\partial D_2$ differs from $\partial D_1$ in the presence of a singular point lying on the boundary of the unbounded component of $\R^2\backslash\overline{(D_1\cup D_2)}$, then the far-field patterns corresponding to $D_j$ and $q_j$ incited by any incoming wave cannot coincide.
\end{theorem}

Theorem \ref{TH2} can be used to distinguish two penetrable scatterers with a piecewise constant potential.
Equivalently, Theorem \ref{TH2} can be reformulated as follows:
\begin{corollary}\label{Coro1}
Let $D_j\subset \R^2$ ($j=1,2$) be two penetrable obstacles in $\R^2$ with the piecewise constant potential functions $q_j$, respectively. Assume that $\partial D_j$ are piecewise analytic and  all non-analytical points of the boundary are singular corners defined by Definition
 \ref{def:S}.  If the far-field patterns corresponding to $(D_j, q_j)$ incited by a single incoming wave are identical, then the boundary of the unbounded component of $\R^2\backslash\overline{(D_1\cup D_2)}$ must be analytic.
 \end{corollary}

Theorem \ref{TH2} and Corollary \ref{Coro1} can be verified in the same manner as the proof of Theorem \ref{TH}. We omit the proofs for simplicity. Note that the above shape identification problem is a formally-determined inverse issue. If the far-field data is available for all incident directions but at fixed energy, uniqueness was verified based on the idea of Isakov; see \cite{Isakov90, Kirsch93}. We also refer to \cite{B2008,IY2012} and \cite[Chapter 10]{CK} for unique determination of potential functions from the data of infinitely many plane waves or the Dirichlet-to-Neumann map.

%

\section{Strongly singular corners always scatter}

This section is devoted to the proof of Theorem \ref{TH} when $\partial D$ contains a strongly singular corner. In the case that $D$ is a polygon, the proof simply follows from \cite{ElHu2015} where variable analytical potential functions were treated.
If $\partial D$ contains a curvilinear corner, the proof was given in \cite{EH2017}. We shall present a proof valid for all strongly singular corners in 2D, under the assumption that $q$ is a piecewise constant potential.

Without loss of generality, we suppose that $O=(0,0)$ is a strongly singular point lying on $\partial D$. Assuming that $u^{sc}$ vanishes in $D^e$, we shall derive a contradiction. Suppose that the boundary $\partial D$ in a neighborhood of $O$ can be expressed as $\Gamma=\{(x_1,f(x_1)): x_1\in(-1,1)\}$, where
$f(x_1)\in C([-1,1])$ is piecewise analytic in $(-1,0]\cup [0,1)$, $f(0)=0$ and $f'(x_1)$ is discontinuous at $x_1=0$.
Since $u^{sc}=0$ in $D^e$,  the Cauchy data of $u$ on $\Gamma$ coincide with those of $u^{in}$, which are analytic. Observing that $q$ is a constant on $\overline{D}$ and $\Gamma$ is piecewise analytic,  by Cauchy-Kovalevskaya theorem, one may extend $u$ analytically from
 $\overline{D}\cap B_1$ to a small neighborhood of $O$ in
 the exterior domain
$D^e\cap B_1$.  For notational convenience, we suppose the extended domain contains $B_1$. Further, the extended function, which we still denote by $u$,  satisfies the Helmholtz equation
\ben
\Delta u+k^2 q_0 u=0\quad\mbox{in}\quad B_1.
\enn
Hence, we deduce  the transmission problem for the Helmholtz equations
\be\label{transmission}\left\{\begin{array}{lll}
\Delta u_j+q_ju_j=0, \quad j=1,2,&\mbox{ in}\quad B_1,\\
u_1=u_2,\quad\frac{\partial u_1}{\partial \nu}= \frac{\partial u_2}{\partial \nu}&\mbox{ on}\quad \Gamma,
\end{array}\right.\en
where
\ben
u_1=u^{in},\quad u_2=u,\qquad q_1=k^2,\quad q_2=k^2 q_0.
\enn
To prove Theorem \ref{TH} for strongly singular corners, it is essential to prove that

\begin{lemma} \label{lem:SS}
  Suppose that $u_j\in H^2(B_1)$ ($j=1,2$) are solutions to (\ref{transmission}). If  $q_1\neq q_2 $, then  $u_1=u_2\equiv0$ in $B_1$.
  \end{lemma}
  By Lemma \ref{lem:SS} and the unique continuation, $u^{in}$ vanishes identically in $\R^2$ which is impossible.
  Hence,  a piecewise constant potential with a strongly singular point on the boundary of the support always scatter.
  The proof of Lemma \ref{lem:SS} follows from an adaption of the arguments in the proof of \cite[Proposition A. 3]{EH2017} to an in homogeneous Helmholtz equation with vanishing Cauchy data.
  To make this paper self-contained,  we present the proof as follows.
\begin{proof} Setting $u:=u_1-u_2$. Then $u$ is a solution to an inhomogeneous Helmholtz equation with vanishing Cauchy data on $\Gamma$:
\be\label{eq:26}\left\{\begin{array}{lll}
\Delta u+q_1u=(q_1-q_2)u_2  &\mbox{in}\quad B_1,\\
u=\frac{\partial u}{\partial \nu}=0&\mbox{on}\quad \Gamma,\\
\Delta u_2+q_2u_2=0 &\mbox{in}\quad B_1.
\end{array}\right.\en
Note that $u$ and $u_2$ are both real-analytic functions in $B_1$.
Denote by $\tilde{\tau}_j(x_1)$ and $\tilde{\nu}_j(x_1)$ ($j=1,2$) the unit tangential and normal vectors on the curves
\[\Gamma_1:=\{(x_1,f(x_1)): x_1\in[0,1)\},\quad\Gamma_2:=\{(x_1,f(x_1): x_1\in(0,-1]\},
\]
which intersect at the corner $O$. Since $f'(x_1)$ is discontinuous at $x_1=0$, the tangential and normal vectors at the corner point, which we denote by
$\tau_j:=\tilde{\tau}_j(0)$ and $\mu_j=\tilde{\nu}_j(0)$,
 are linearly independent.  Without loss of generality we suppose that $\nu_1=a_1\tau_1+a_2\tau_2$ with $a_1,a_2\in\R$, $a_2\neq 0$. Hence,
\be\label{Derivative}
\partial_{\tau_2}=\frac{1}{a_2}\partial_{\nu_1}-\frac{a_1}{a_2}\partial_{\tau_1}.
\en
We shall prove by induction  that $\nabla^m u(O)=0$ for all $m\in\N_0$, which implies the lemma.

From the Dirichlet and Neumann boundary conditions of $u$ on $\Gamma=\Gamma_1\cup \Gamma_2$ we see that
\be\label{D-12}
u=\nabla u=0,\quad
\partial_{\tau_1}^2u=\partial_{\tau_2}^2u=\partial_{\nu_1}\partial_{\tau_1}u=0\quad\mbox{at the corner}\; O.
\en
Combining (\ref{Derivative}) and (\ref{D-12}) gives the relation
$\partial_{\tau_1}\partial_{\tau_2}u=0$ at $O$.
Since each entry of the vector $\nabla^2$ can be expanded as a linear combination of $\partial_{\tau_1}^2$, $\partial^2_{\tau_2}$ and $\partial_{\tau_1}\partial_{\tau_2}$, we obtain
\be\label{eq:27}
\nabla^2u=0\quad \mbox{at}\quad O.
\en
Consequently, it follows from the equations in (\ref{eq:26}) that
\be\label{eq:28}
u_2=\Delta u_2=0\quad \mbox{at}\quad O,
\en
where we have used the assumption that $q_1\neq q_2$.

To prove that $\nabla^3u(O)=0$, we observe that
\ben
\partial_{\tau_1}^3u=\partial_{\tau_2}^3u=
\partial_{\tau_1}^2\partial_{\nu_1}u=
\partial_{\tau_2}^2\partial_{\nu_2}u=0 \quad\mbox{at}\; O.
\enn
Applying $\partial_{\tau_1}^2$ to both sides of (\ref{D-12}) yields $\partial_{\tau_1}^2\partial_{\tau_2}u(O)=0$. Analogously we can get $\partial_{\tau_2}^2\partial_{\tau_1}u(O)=0$.
Hence, the relation $\nabla^3u(O)=0$ follows from the fact that the differential operators $\partial_{\tau_1}^3, \partial_{\tau_1}^2\partial_{\tau_2}, \partial_{\tau_1}\partial_{\tau_2}^2$ and $\partial_{\tau_2}^3$ span the vector $\nabla^3$. Taking $\nabla$ on the equations in (\ref{eq:26}) gives
\be\label{eq:29}
\nabla u_2=\nabla\Delta u_2=0\quad \mbox{at}\quad O.
\en

Now we want to verify that $\nabla^4u(O)=0$. Arguing as in the previous step we get
\be\label{eq1:8}
\partial_{\tau_1}^4u=\partial_{\tau_1}^3\partial_{\tau_2}u=\partial_{\tau_1}\partial_{\tau_2}^3u=
\partial_{\tau_2}^4u=0\quad\mbox{at}\; O.
\en Hence it suffices to prove $\partial_{\tau_1}^2\partial_{\tau_2}^2u(O)=0$.
Taking $\Delta$ on the first equation in (\ref{eq:26}) and using (\ref{eq:27})-(\ref{eq:28}), we find
\ben
\Delta^2 u(O)=(q_1-q_2)\Delta u_2(O)-q_1 \Delta u(O)=0.
\enn
On the other hand, using (\ref{eq1:8}) and $\partial_{\nu_1}=a_1\partial_{\tau_1}+a_2\partial_{\tau_2}$, we deduce that
\ben
\Delta^2u(O)=[\partial_{\nu_1}^2 +\partial_{\tau_1}^2]^2u(O)=[2(1+a_1^2)a_2^2+4a_1^2a_2^2]\,\partial_{\tau_1}^2\partial_{\tau_2}^2u(O),
\enn from which the relation $\partial_{\tau_1}^2\partial_{\tau_2}^2u(O)=0$ follows. This proves $\nabla^4u(O)=0$.
Now, differentiating the equations in (\ref{eq:26}) yields
\ben
\nabla^2 u_2=\nabla^2\Delta u_2=0\quad \nabla\Delta^2 u=0\quad \mbox{at}\quad O.
\enn

For $m>4$, we make the induction hypothesis that
\be\label{induction-hypothesis}
\nabla^j u(O)=\nabla^{j-3}\Delta^2 u(O)=0,\quad  \nabla^{j-2}u_2(O)=\nabla^{j-2}\Delta u_2(O)=0\quad\mbox{for all}\;j=0,1,\cdots,m.
\en We then need to verify that the above relations hold for $j=m+1$, that is,
\ben
\nabla^{m+1} u(O)=\nabla^{m-2}\Delta^2 u(O)=0,\quad  \nabla^{m-1}u_2(O)=\nabla^{m-1}\Delta u_2(O)=0.
\enn

We first prove $\nabla^{m+1}u =0$ at $O$.
For $j\in \N_0$, denote by $\nabla^j_\tau$ the vector of all tangential derivatives of order $j$, i.e.,
\ben
\nabla^j_\tau u=\left\{\partial_{\tau_1}^{j_1}\partial^{j_1}_{\tau_2}\,u: \quad j_1,j_2\in\N_0,j_1+j_2=j\right\}.
\enn In view of the vanishing of the Cauchy data on $\Gamma$ and using (\ref{Derivative}) again, we have
\ben
\nabla_\tau^{m-3}\Delta^2u=\partial_{\tau_1}^{m+1}u=\partial_{\tau_1}^m\partial_{\tau_2}u
=\partial_{\tau_1}\partial_{\tau_2}^mu=\partial_{\tau_2}^{m+1}u=0\quad\mbox{at}\quad O.
\enn
It was shown in \cite[Proposition A.3]{EH2017} that the span of the differential operators $\nabla_\tau^{m-3}\Delta^2$,
$\partial_{\tau_1}^{m+1}$, $\partial_{\tau_1}^m\partial_{\tau_2}$,
$\partial_{\tau_1}\partial_{\tau_2}^m$ and $\partial_{\tau_2}^{m+1}$ contains the vector $\nabla^{m+1}_\tau$. Hence, the relation $\nabla^{m+1}u =0$ at $O$ follows.
Taking $\nabla^{m-1}$ on the equations in (\ref{eq:26}) and using $q_1\neq q_2$, we see
\ben
\nabla^{m-1} u_2=\nabla^{m-1}\Delta^2 u_2=0\quad \mbox{at}\quad O.
\enn
 Taking $\nabla^{m-1}\Delta^2$ on the first equation in (\ref{eq:26}) gives $\nabla^{m-2}\Delta^2 u(O)=0$.
 By induction we obtain $u_1=u_2\equiv 0$ in $B_1$.
\end{proof}

\section{Weakly singular corners always scatter}
To prove Theorem \ref{TH} for weakly singular points lying on $\partial D$, we only need to show that

\begin{lemma} \label{lem:WS}
Let  $\Gamma$ be the profile of the function (\ref{f1}) .
  Suppose that $u_j$ ($j=1,2$) are solutions to the Helmholtz equation $\Delta u_j+q_ju_j=0, j=1,2 $ in $B_1$ with  $q_1\neq q_2 $, subject to the transmission conditions  $ u_1=u_2,\frac{\partial u_1}{\partial \nu}= \frac{\partial u_2}{\partial \nu}$ on $\Gamma$. If  $O\in \Gamma$ is a weakly singular point and $q_1\neq q_2$, then  $u_1=u_2\equiv0$.
  \end{lemma}
  It seems non-trivial to prove the above lemma by extending the analysis in the proof of Lemma \ref{lem:SS} to the case of weakly singular corners. The analytical approach of using polar coordinates (see \cite{ElHu2015}) also turns out to be complicated. Below we shall present a novel approach by using the expansion of solutions to the Helmholtz equation in the Cartesian coordinate system.

Since  $u_j$ satisfies the Helmholtz equation and $q_j$  is constant, the solution $u_j$ is analytic in $B_1$. Hence,
$u_j$ can be expanded into the convergent Taylor expansion
\ben
 u_j=\sum_{n\ge0} \sum_{m \ge0} a^{(j)}_{n,m}\, x_1^n\, x_2^m \quad\mbox {in} \quad B_1,
 \enn
 where the coefficients $a^{(j)}_{n,m}\in \C$ satisfy the relation
\begin{equation}
(n+1)(n+2) a^{(j)}_{n+2,m}+(m+1)(m+2)a^{(j)}_{n,m+2}+q_j a^{(j)}_{n,m}=0.
\label{eq:1}
\end{equation}
\par
Set $u=u_1-u_2$. Then $u$ admits the Taylor expansion
$$u=\sum_{n\ge0} \sum_{m \ge0} a_{n,m} x_1^n x_2^m\quad\mbox{in}\quad B_1,\qquad a_{n,m}:=a_{n,m}^{(1)}-a_{n,m}^{(2)},$$
and satisfies the equation $\Delta u+ q_1 u=(q_2-q_1)u_2$ in $B_1$. The implies that  the coefficients $a_{n,m}$ fulfills the recursive relation
\begin{equation}
(n+1)(n+2) a_{n+2,m}+(m+1)(m+2)a_{n,m+2}+q_1 a_{n,m}=(q_2-q_1)a^{(2)}_{n,m}.
\label{eq:2}
\end{equation}
Combining (\ref{eq:1}) and (\ref {eq:2}), we deduce that
\begin {equation}
\begin{split}
0&=(m+4)(m+3)(m+2)(m+1)a_{n,m+4}+(n+4)(n+3)(n+2)(n+1)a_{n+4,m}\\
&\quad+2(n+2)(n+1)(m+2)(m+1)a_{n+2,m+2}\\
&\quad+(q_1+q_2)(n+2)(n+1)a_{n+2,m}
+(q_1+q_2)(m+2)(m+1)a_{n,m+2}\\
&\quad+q_2q_1 a_{n,m}.
\end{split}
\label{eq:3}
\end{equation}
We shall prove $a_{n,m}=0$ for all $n, m\in \N$ through (\ref{eq:3}) and the transmission conditions
\ben
u=\partial_\nu u=0\quad\mbox{on}\quad \Gamma.
\enn
This together with  (\ref{eq:2}) would give rise to $a^{(2)}_{n,m}=a^{(1)}_{n,m}=0.$

Denote by $\Gamma_1:=\{(x_1, f(x_1): x_1\in[0,1)\}$ and $\Gamma_2=\{(x_1, f(x_1): x_1\in(-1,0])\}$, with the normal directions given by
\[
\nu^{(1)}(x_1)=(\alpha_1 c_1x_1^{\alpha_1-1},-1)^\top, \quad x_1>0;\qquad
\nu^{(2)}(x_1)=(\alpha_2 c_2x_1^{\alpha_2-1},-1)^\top,\quad x_1<0,
\] respectively. Observe that
\[
\frac{\partial u}{\partial x_1}=\sum_{n\ge0}\sum_{m\ge0}a_{n+1,m}(n+1)x_1^n x_2^m,\quad
\frac{\partial u}{\partial x_2}=\sum_{n\ge0}\sum_{m\ge0}a_{n,m+1}(m+1)x_1^n x_2^m.
\]
It follows from $\partial_\nu u=0$ on $\Gamma_j$ ($j=1,2$) that
\be\label{eq:4}
 \alpha_1\sum_{n\ge0}\sum_{m\ge0}a_{n+1,m}(n+1)x_1^{n+\alpha_1 m+\alpha_1-1} c_1^{m+1}
-\sum_{n\ge0}\sum_{m\ge0}a_{n,m+1}(m+1)x_1^{n+\alpha_1 m}c_1^{m}=0, \\ \label{eq:5}
\alpha_2 \sum_{n\ge0}\sum_{m\ge0}a_{n+1,m}(n+1)x_1^{n+\alpha_2 m+\alpha_2-1} c_2^{m+1}
-\sum_{n\ge0}\sum_{m\ge0}a_{n,m+1}(m+1)x_1^{n+\alpha_2 m}c_2^{m}=0.
\en

Without loss of generality, we suppose that $\alpha_1\leq \alpha_2$. In order to prove Lemma \ref{lem:WS}, we will consider two cases:
\par
\indent Case 1: $ \alpha_1=\alpha_2\geq 2$;\qquad\qquad\qquad\qquad\qquad
Case 2: $2\leq\alpha_1< \alpha_2$.

The proofs in  Cases 1 and 2
will be carried out in the subsequent two subsections, separately.

\subsection{Proof of Lemma \ref{lem:WS} when $\alpha_1=\alpha_2\geq 2$.}\label{subsec}

For notational convenience we set $\alpha:=\alpha_1=\alpha_2\geq 2$.
Equating coefficients of ${x_1}^l$ ($l\in \N, l\ge \alpha-1$) in (\ref{eq:4}) and (\ref{eq:5}) and changing properly the summation indices,   we  obtain
\be\label{eq:6}
\sum_{n+\alpha m=l-\alpha+2,n\ge1, m\ge0}\alpha n a_{n,m}c_1^{m+1}-\sum_{n+\alpha m=l+\alpha,n\ge0, m\ge1}ma_{n,m}c_1^{m-1}=0,\\ \label{eq:7}
\sum_{n+\alpha m=l-\alpha +2,n\ge1, m\ge0}\alpha n a_{n,m}c_2^{m+1}-\sum_{n+\alpha m=l+\alpha ,n\ge0, m\ge1}ma_{n,m}c_2^{m-1}=0.
\en On the other hand, the Dirichlet condition
$u=0$ on $\Gamma$ gives  the relations
\begin{equation}
 \sum_{n\ge0}\sum_{m\ge0}a_{n,m}x^{n+\alpha m}_1{c_1}^{m}=0,\quad
 \sum_{n\ge0}\sum_{m\ge0}a_{n,m}x^{n+\alpha m}_1c_2^m=0,
\label{eq:8}
\end{equation}
Equating coefficients of ${x_1}^l$ ($l\in \N$, $l\geq 0$) in (\ref{eq:8}), we get
\be\label{eq:9}
\sum_{n+\alpha m=l,n\ge0,m\ge0}a_{n,m}{c_1}^m=0,\quad
\sum_{n+\alpha m=l,n\ge0,m\ge0}a_{n,m}{c_2}^m=0.
\en
By (\ref{eq:6}),(\ref{eq:7}) and ({\ref{eq:9}}) we shall prove $a_{n,m}=0$ for all $ n+\alpha m=j$ ($j\in \N$)  by an induction argument on the index $j\in \N$. We divide the proof into four steps.
\par
{\bf Step 1:}  Prove $a_{j,0}=0$ for all $j=0,1,\ldots \alpha-1$. This follows from
(\ref{eq:9}) with $n+\alpha m=l$ for $l=0,1,\ldots \alpha-1$.
\par
{\bf Step 2:} Prove $a_{j,0}=a_{j-\alpha,1}=0$ for all $j=\alpha, \alpha+1,\cdots, 2\alpha-1$.
Setting $ l=\alpha,\alpha+1,\dots,2\alpha-1$ in (\ref{eq:9}),  we obtain
\[
a_{l,0}+c_1a_{l-\alpha,1}=0,\quad
a_{l,0}+c_2a_{l-\alpha,1}=0.
\]
Since $c_1\ne c_2$, we see $a_{l,0}=a_{l-\alpha,1}=0$.

{\bf Step 3:}  Prove $a_{j,0}=a_{j-\alpha,1}=a_{j-2\alpha, 2}=0$ for all $j=2\alpha, 2\alpha+1,\cdots, 3\alpha-1$.
As done in previous two steps, setting $ n+\alpha m=2\alpha, 2\alpha+1,\dots,3\alpha-1$ in (\ref{eq:9}), we find
\be\label{neq:1}
a_{j,0}+c_1a_{j-\alpha,1}+{c_1}^2\,a_{j-2\alpha,2}=0,\quad
a_{j,0}+c_2a_{j-\alpha,1}+{c_2}^2\,a_{j-2\alpha,2}=0.
\en
On the other hand, one may conclude from Steps 1 and 2 that
\[
a_{n,m}=0   \quad \mbox{if}\quad \ n+\alpha m < j.
\]
This together with $\alpha\geq 2$ implies that
\ben
0=\sum_{n+\alpha m=j-2\alpha+2,n\ge1, m\ge0}\alpha n a_{n,m}c_1^{m+1}=\sum_{n+\alpha m=j,n\ge1, m\ge0}\alpha n a_{n,m}c_2^{m+1}
\enn
Hence, setting $l=j-\alpha$ in (\ref{eq:6}) and (\ref{eq:7}) gives the relations
\be\label{neq:2}
a_{j-\alpha,1}+2{c_1}a_{j-2\alpha,2}=0,\quad
a_{j-\alpha,1}+2{c_2}a_{j-2\alpha,2}=0.
\en
Therefore,  combining (\ref{neq:1}) and (\ref{neq:2})  yields $a_{j,0}=a_{j-\alpha,1}=a_{j-2\alpha,2}=0$. Further, we conclude from Steps 1-3 that
\be\label{neq:3}
a_{n,m}=0   \quad \mbox{if}\quad \ n+\alpha m < 3\alpha.
\en

{\bf Step 4:} Prove $a_{j,0}=a_{j-\alpha,1}=a_{j-2\alpha, 2}=a_{j-3\alpha, 3}=0$ for all $j=3\alpha, 3\alpha+1,\cdots, 4\alpha-1$.

Setting $ n+\alpha m=3\alpha,3\alpha+1,\dots,4\alpha-1$ in (\ref{eq:9}),  we get for such $j$ that
\ben
a_{j,0}+c_1a_{j-\alpha,1}+{c_1}^2a_{j-2\alpha,2}+{c_1}^3a_{j-3\alpha,3}=0,\\
a_{j,0}+c_2a_{j-\alpha,1}+{c_2}^2a_{j-2\alpha,2}+{c_2}^3a_{j-3\alpha,3}=0,
\enn
Setting $l=j-\alpha$ in (\ref{eq:6})-(\ref{eq:7}) and making use of  (\ref{neq:3}), we  obtain
\ben
a_{j-\alpha,1}+2c_1a_{j-2\alpha,2}+3{c_1}^2a_{j-3\alpha,3}=0,\\
a_{j-\alpha,1}+2c_2a_{j-2\alpha,2}+3{c_2}^2a_{j-3\alpha,3}=0.
\enn
For fixed $j\in\{3\alpha, 3\alpha+1,\cdots, 4\alpha-1 \}$, the previous relations can be written as the system
\be\label{neq:4}
\begin {pmatrix}
1&c_1&{c_1}^2&{c_1}^3\\
1&c_2&{c_2}^2&{c_2}^3\\
0&1&2c_1&3{c_1}^2\\
0&1&2c_2&3{c_2}^2
\end {pmatrix} \begin{pmatrix}
a_{j,0} \\ a_{j-\alpha, 1} \\ a_{j-2\alpha, 2} \\ a_{j-3\alpha, 3}
\end{pmatrix}=0.
\en
It is not difficult to check that the determinant of the matrix on the left hand side of (\ref{neq:4})  is $-(c_1-c_2)^4\neq0$, implying that $ a_{j,0}=a_{j-\alpha,1}=a_{j-2\alpha,2}=a_{j-3\alpha,3}=0$.
Hence, it holds that
\ben
a_{n,m}=0   \quad \mbox{if}\quad \ n+\alpha m <4\alpha.
\enn
{\bf Step 5:} Induction arguments.  We make the induction hypothesis
$a_{n,m}=0$ for all $n+\alpha m<M$ for some $M\ge4\alpha, M\in \N$.
We need to prove that
\begin{equation}
a_{n,m}=0\quad \mbox{if}\quad n+\alpha m=M.
\label{eq:11}
\end{equation}
We first claim that
\be\label{neq:5}
a_{n,m}=0, \quad \mbox{if}\quad  n+\alpha m=M, m\ge4.
\en
Let $n'=n, m'=m-4\ge 0$. Then we see  i
\[
n'+\alpha (m'+4)=M, \quad n', m'\geq 0.
\]
One can readily prove that
\[
n'+4+\alpha m'<M,\quad
n'+2+\alpha (m'+2)<M,\quad
n'+2+\alpha m'<M,\quad
n'+\alpha( m'+2)<M.
\]
Therefore, by induction hypothesis,
\[
a_{n'+4,m'}=a_{n'+2,m'+2}=a_{n'+2,m'}=a_{n',m'+2}=0.
\]
Using (\ref{eq:3}), we get the relation
\[
a_{n',m'+4}=0,\ \mbox{if}\quad n'+\alpha(m'+4)=M, \quad n',m' \ge 0,
\] which proves (\ref{neq:5}).

To proceed with the proof we set $ l=M$ in (\ref{eq:9}) to  obtain
\ben
a_{M,0}+c_1a_{M-\alpha,1}+{c_1}^2a_{M-2\alpha,2}+{c_1}^3a_{M-3\alpha,3}=0,\\
a_{M,0}+c_2a_{M-\alpha,1}+{c_2}^2a_{M-2\alpha,2}+{c_2}^3a_{M-3\alpha,3}=0,
\enn where the relation (\ref{neq:5}) was again used. On the other hand, setting $l=M-\alpha$ in (\ref{eq:6})-(\ref{eq:7}) and recalling  the induction hypothesis, we see
\ben
a_{M-\alpha,1}+2c_1a_{M-2\alpha,2}+3{c_1}^2a_{M-3\alpha,3}=0,\\
a_{M-\alpha,1}+2c_2a_{M-2\alpha,2}+3{c_2}^2a_{M-3\alpha,3}=0.
\enn
Note that the coefficient matrix for the unknowns $a_{M,0}, a_{M-\alpha, 1}, a_{M-2\alpha, 2}$ and $a_{M-3\alpha, 3}$  is the same as the 4-by-4 matrix on the left hand side of (\ref{neq:4}). Since the determinant of this matrix does not vanish, we obtain
 $ a_{M,0}=a_{M-\alpha,1}=a_{M-2\alpha,2}=a_{M-3\alpha,3}=0$.
 This together with (\ref{neq:5}) proves (\ref{eq:11}).

 By induction, it holds that $a_{n,m}=0$ for all $n, m\in \N.$ In view of (\ref{eq:2}) and the condition $q_1\neq q_2$, we obtain $a_{n,m}^{(1)}=a_{n,m}^{(2)}=0$ for all $n, m\in \N$. Finally, we get $u_1=u_2\equiv0$ in $B_1$ by the analyticity of $u_j$ ($j=1,2$).

\subsection{Proof of Lemma \ref{lem:WS} when $2\leq\alpha_1<\alpha_2$.}

We first observe that the powers of $x_1$ in the first summation on the left hand side of
(\ref{eq:4}) and (\ref{eq:5}) are all greater than or equal to $\alpha_j-1$, whereas those in the second summation start from zero. Hence,
equating coefficients of the term $x_1^{l}$ ($l<\alpha_j-1$) in (\ref{eq:4}) and (\ref{eq:5}) yields
\be\label{eq:12}
\sum_{n\ge0,m\ge1,n+\alpha_1m=l+\alpha_1}a_{n,m}mc_1^{m-1}=0,\quad l\le\alpha_1-2,\\ \label{eq:14}
\sum_{n\ge0,m\ge1,n+\alpha_2m=l+\alpha_2}a_{n,m}mc_2^{m-1}=0,\quad l\le\alpha_2-2.
\en
Analogously, equating coefficients of the term $x_1^{l}$ for $l\ge \alpha_j-1$, we obtain
\be\label{eq:13}
\sum_{n+\alpha_1m=l-\alpha_1+2,n\ge1, m\ge0}\alpha_1 na_{n,m}c_1^{m+1}-\sum_{n+\alpha_1m=l+\alpha_1,n\ge0, m\ge1}ma_{n,m}c_1^{m-1}=0, l\ge \alpha_1-1,\\ \label{eq:15}
\sum_{n+\alpha_2m=l-\alpha_2+2,n\ge1, m\ge0}\alpha_2 na_{n,m}c_2^{m+1}-\sum_{n+\alpha_2m=l+\alpha_2,n\ge0, m\ge1}ma_{n,m}c_2^{m-1}=0,l\ge \alpha_2-1.
\en
From the Dirichlet boundary condition $u=0$ on $\Gamma$, we obtain
\[
\sum_{n\ge0}\sum_{m\ge0}a_{n,m}x^{n+\alpha_1 m}_1{c_1}^{m}=0,\quad \sum_{n\ge0}\sum_{m\ge0}a_{n,m}x^{n+\alpha_2 m}_1c_2^m=0,
\]
which implies that
\be\label{eq:16}
\sum_{n+\alpha_1m=l,n\ge0,m\ge0}a_{n,m}{c_1}^m=0,\quad l\ge0, \\\label{eq:17}
\sum_{n+\alpha_2m=l,n\ge0,m\ge0}a_{n,m}{c_2}^m=0, \quad l\ge0.
\en
Since $\alpha_1<\alpha_2$, the proof in this section is more complicated than previous subsection. We shall still apply the induction argument to prove that $a_{n,m}=0$ for all $ n+\alpha m=j$, $j\in \N$. Below we carry out the proof under the assumption that $c_1\neq 0$. If $c_1=0$, we have $c_2\neq 0$ by assumption. Then the interface can be locally parameterized by the function given in (\ref{f1}) with $\alpha_1=\alpha_2$. Hence, the vanishing of $u_j$ in $B_1$ follows from the same arguments used in subsection \ref{subsec}, where the case $c_1=0$  is covered.

\textbf{Step 1:} Prove $a_{j,0}=0$ for all $j=\alpha_1,\alpha_1+1,\ldots, \alpha_2-1 $. This follows from (\ref{eq:17}).

\textbf{Step 2:} Prove $a_{j,0}=a_{j-\alpha_1,1}=0$ when $j=\alpha_1, \alpha_2, \ldots, min ({\alpha_2-1,2\alpha_1-1})$.

Setting $j=\alpha_1,\ldots 2\alpha_1-1$ in (\ref{eq:16}), we obtain
\begin{equation}
a_{j,0}+a_{j-\alpha_1,1}c_1=0.
\label {eq:18}
\end{equation}
This together with the condition $c_1\neq 0$ and the fact that $a_{l,0}=0$ for $l=\alpha_1,\alpha_1+1,\ldots, \alpha_2-1 $ (see Step 1) gives the desired result.

\textbf{Step 3:} Induction arguments. Assuming that
 \ben
 a_{n,m}=0\quad\mbox{for all}\quad n+m\alpha_1<M,\quad M> min({\alpha_2-1,2\alpha_1-1}),
 \enn
we will prove that
\begin{equation}
a_{n,m}=0\quad\mbox{for all}\quad n+m\alpha_1=M.
\label{eq:19}
\end{equation}
Setting $l=M$ in (\ref{eq:17}) gives
\begin{equation}
\sum_{n+m \alpha_2=M,n\ge0,m\ge0}a_{n,m}c_2^m=0.
\label{eq:20}
\end{equation}
Since $\alpha_1<\alpha_2$, the indices $n, m$ appearing in the above summation fulfill
$n+m\alpha_1<M$ if $m\ne0$.
By induction hypothesis, this implies that
\[
a_{n,m}=0, \quad \text{if} \ n+m\alpha_2=M, m\ne0.
\]
When $m=0$, it follows from (\ref{eq:20}) that $a_{M,0}=0$. Now, it remains to prove
\be\label{neq:6}
a_{n,m}=0 \quad \text{for all}\ n+\alpha_1m=M, \quad m\ne0.
\en
in the following cases.
\subsubsection{Case 1: $M\le2\alpha_1-1.$}
Setting $l=M$ in (\ref{eq:16}), we obtain
\[
a_{M,0}+a_{M-\alpha_1,1}c_1=0,
\]
which together with $ a_{M,0}=0$ and $c_1\neq 0$ leads to  $a_{M-\alpha_1,1}=0$.
This proves (\ref{neq:6}) when $M\le2\alpha_1-1.$

\subsubsection{Case 2: $2\alpha_1\le M\le3\alpha_1-1$.}
 Letting $l=M$ in (\ref{eq:16}) and using again the fact that $a_{M,0}=0$, we obtain
\be\label{neq:7}
a_{M-\alpha_1,1}c_1+a_{M-2\alpha_1,2}{c_1}^2=0.
\en
Setting $l=M-\alpha_1$ in (\ref{eq:13}) and making use of the induction hypothesis
\[
a_{n,m}=0,\quad\text{for all}\; n+\alpha_1m=l-\alpha_1+2<M,
\]
we obtain
\be\label{neq:8}
a_{M-\alpha_1,1}+2a_{M-2\alpha_1,2}{c_1}=0.
\en
Combining (\ref{neq:7}) and (\ref{neq:8}) leads to $a_{M-\alpha_1,1}=a_{M-2\alpha_1,2}=0$,
which proves (\ref{neq:6}).

\subsubsection{Case 3: $3\alpha_1\le M \le 4\alpha_1-1 $.}\label{sub-1}
 As done in previous cases, setting $l=M$ in (\ref{eq:16}) and using $a_{M,0}=0$ gives
\begin {equation}
a_{M-\alpha_1,1}c_1+a_{M-2\alpha_1,2}{c_1}^2+a_{M-3\alpha_1,3}{c_1}^3=0.
\label{eq:21}
\end{equation}

Setting $l=M-\alpha_1$ in (\ref{eq:13}). Recalling from the induction hypothesis that
\[
a_{n,m}=0,\quad\text{if}\ n+\alpha_1m=l-\alpha_1+2<M,
\]
we obtain
\begin{equation}
a_{M-\alpha_1,1}+2a_{M-2\alpha_1,2}{c_1}+3 a_{M-3 \alpha_1,3}{c_1}^2=0.
\label{eq:22}
\end{equation}

Next we will show $a_{M-\alpha_1,1}=0$. Write $N=M-\alpha_1+\alpha_2$ for notational simplicity.

 If $N\le2\alpha_2-2$, setting $l=N-\alpha_2$ in (\ref{eq:14}) gives
$a_{M-\alpha_1,1}=0.$

If $N\ge2\alpha_2-1$,  we have the relation
\[
M-\alpha_1+\alpha_2-2\alpha_2+2=M-\alpha_1-\alpha_2+2 <M.
\]
Letting $l=M-\alpha_1$ in (\ref{eq:15}), we can obtain
\begin{equation}
\sum_{n+\alpha_2m=N,n\ge0, m\ge1}a_{n,m}m{c_2}^{m-1}=0,
\label{eq:23}
\end{equation}
Since $\alpha_2>\alpha_1$, it holds that
\[
 n+\alpha_1m <M\quad\mbox{for all}\quad n+\alpha_2m=N, m\ge2,
\]
implying that
\[
a_{n,m}=0\quad\mbox{for all}\quad n+\alpha_2m=N, m\ge2.
\]
due to the induction hypothesis.
Hence, it follows from (\ref{eq:23}) that $a_{M-\alpha_1,1}=0.$
Combining this with (\ref{eq:21}) and (\ref{eq:22}) and the fact that $c_1\ne0$, we obtain  $a_{M-2\alpha_1,2}=a_{M-3\alpha_1,3}=0$, which proves (\ref{neq:6}).

\subsubsection{Case 4:  $M\ge 4\alpha_1$.}\label{sec:4.2.4}
 We first prove that
\begin{equation}
a_{n,m+4}=0, \quad\text{if} \ n+\alpha_1(m+4)=M, n\ge0, m\ge0.
\label{eq:24}
\end{equation}

Supposing that indices $n,m\ge 0$ in (\ref{eq:3}) satisfy the relation $n+\alpha_1 (m+4)=M$.
 Then we have
\ben
n+4+\alpha_1 m<M,\quad
n+2+\alpha_1 (m+2)<M, \\
n+2+\alpha_1 m<M,\quad
n+\alpha_1( m+2)<M,
\enn
By induction hypothesis, we see
\[
a_{n+4,m}=a_{n+2,m+2}=a_{n+2,m}=a_{n,m+2}=0.
\]
Hence, the relation (\ref{eq:24}) follows from (\ref{eq:3}).
To prove (\ref{neq:6}) we only need to verify
 \begin{equation}
 a_{M-\alpha_1,1}=a_{M-2\alpha_1,2}=a_{M-3\alpha_1,3}=0.
 \label{eq:25}
 \end{equation}
 \par
 Analogously to the Case 3, we will show that $a_{M-\alpha_1,1}=0$ by setting $N=M-\alpha_1+\alpha_2$.

If $N\le2\alpha_2-2$, letting $l=N$ in (\ref{eq:14}) leads to $a_{M-\alpha_1,1}=0$.

If $N\ge2\alpha_2-1$, letting $l=M-\alpha_1$ in (\ref{eq:15}) and noting that
\[
 M-\alpha_1+\alpha_2-2\alpha_2+2=M-\alpha_1-\alpha_2+2 <M,
\]
we obtain
\be\label{neq:9}
\sum_{n+\alpha_2m=N,n\ge0, m\ge1}a_{n,m}m{c_2}^{m-1}=0.
\en
Similar to the arguments in subsection \ref{sub-1}, we can obtain using the  induction hypothesis that
\[
a_{n,m}=0\quad \mbox{for all}\quad n+\alpha_2m=N, m\ge2,
\]
because $n+\alpha_1m<M$ for such indices $n$ and $m$.
Therefore,  we get $a_{M-\alpha_1,1}=0$ from (\ref{neq:9}).
Now, setting $l=M$  in (\ref{eq:16}), using (\ref{eq:24}) and the fact that $a_{M,0}=a_{M-\alpha_1,1}=0$, we see
\be\label{neq:11}
a_{M-2\alpha_1,2}{c_1}^2+a_{M-3\alpha_1,3}{c_1}^3=0.
\en
On the other hand, setting $l=M-\alpha_1$ in (\ref{eq:13}),  using  $a_{M-\alpha_1,1}=0$ and the relations
\[
a_{n,m}=0\quad\mbox{for all}\quad n+\alpha_1m=l-\alpha_1+2=M-2\alpha_1-2<M,
\]
we deduce that
\be\label{neq:10}
2a_{M-2\alpha_1,2}{c_1}+3 a_{M-3 \alpha_1,3}{c_1}^2=0.
\en
Since $c_1\ne0$, we obtain $a_{M-2\alpha_1,2}=a_{M-3\alpha_1,3}=0$ by combining (\ref{neq:11}) and (\ref{neq:10}). This finishes the proof of (\ref{neq:6}) when $M\geq 4\alpha_1$.

Finally, the relation (\ref{eq:19}) follows from (\ref{neq:6}) and the fact that
 $a_{M,0}=0$. The proof of Lemma \ref{lem:WS} is thus complete under the assumption that $2\leq\alpha_1<\alpha_2$.

\section{Concluding remarks}

We remark that Lemma \ref{lem:WS} does not hold true if  the curve $\Gamma$ is analytic. Counterexamples can
be easily constructed when $\Gamma$ is a line segment (see  \cite[Remark 3.3]{EH2017} ) or a circle.
If $\Gamma\subset B_1$ is a circle of radius $R<1$ centered at the origin,  one may find interior transmission eigenvalues (ITEs) ( or equivalently,  $q_1$ and $q_2$) such that the coupling problem
\be\label{ITE}
\Delta u_j + q_j u_j =0\quad\mbox{in}\quad \Omega,\quad u_1=u_2,\quad \partial_\nu u_1=\partial_\nu u_2\quad\mbox{on}\quad \Gamma,
\en
admits non-trivial solutions $u_1$ and $u_2$ in $B_R$ (see e.g.,\cite{CPJ}), which can be analytically extended to $B_1$. Here $\Omega\subset B_1$ denotes the domain enclosed by the closed curve $\Gamma$. Our Lemmas \ref{lem:SS} and \ref{lem:WS} imply that, if $\Gamma$ possesses a singular point,  the non-trivial solutions $u_j$ to (\ref{ITE}) can not be analytically extended onto $\overline{B}_1$.
We refer to \cite{CGH2010, CP, CPJ, PJ,S} for the existence of ITEs in inverse scattering theory. Note that all results of this paper carry over to variable potential functions which is a constant in a small neighborhood of the
singular point under question. The singular points considered here form only a subset of non-analytical points of $\Gamma$.
 We conjecture that Lemmas \ref{lem:SS} and \ref{lem:WS} remain valid under the weak assumption that $\Gamma$ contains a single non-analytical point. However, the proof requires novel mathematical arguments and the progress along this direction will be reported in our forthcoming publications.

\section{Acknowledgments}
G. Hu would like to thank J. Elschner for private discussions on the proof of Lemma \ref{lem:SS}.
The work of G. Hu is supported by the NSFC grant (No. 11671028), NSAF grant (No. U1530401) and the 1000-Talent Program of Young Scientists in China.
The work of L. Li and J. Yang are partially supported by the National Science Foundation of
China (61421062, 61520106004) and Microsoft Research of Asia.

\end{document}